\documentclass[12pt,a4paper]{article}

\usepackage{amsmath}
\usepackage{amsfonts, amssymb}
\usepackage{epsfig}

\def\R{\mathbb{R}}
\def\E{\mathbb{E}}

\def\tr{{\tilde R}}
\def\b{{\mathbf B}}

\def\K{{\textsc{k}}}
\def\kh{{\widehat{\K}}}
\def\H{{\textsc{h}}}
\def\sgn{{\textrm{sgn}\,}}
\def\W{{\mathcal W}}

\newtheorem{thm}{Theorem}[section]

\newtheorem{prop}[thm]{Proposition}

\newtheorem{rem}[thm]{Remark}

\newcommand{\pf}{\noindent {\it Proof:\ }}
\newcommand{\remark}{\noindent {\bf Remark:\ }}

\def\be#1\ee{\begin{equation}#1\end{equation}}

\begin{document}

\title{\bf Bifractional Brownian Motion: \\
Existence and Border Cases}
\author{
   Mikhail Lifshits
   \footnote{St.Petersburg State University, Russia, Stary
   Peterhof, Bibliotechnaya pl.,2,
   email {\tt mikhail@lifshits.org} and Dept. Math., Link\"oping University.}
   \thanks{The work of both authors was supported by grants NSh.2504.2014.1, RFBR 13-01-00172,
and SPbSU 6.38.672.2013.}
   \and
    Ksenia Volkova
   \footnote{St.Petersburg State University, Russia, Stary
   Peterhof, Bibliotechnaya pl.,2,
   email \ {\tt efrksenia@gmail.com}}
   }
\date{\today}
\maketitle

\begin{abstract} Bifractional Brownian motion (bfBm) is a centered Gaussian
process with covariance
\[
   R^{(\H,\K)}(s,t)=
   2^{-\K} \left( \left(|s|^{2\H}+|t|^{2\H} \right)^{\K}-|t-s|^{2\H\K}\right),
   \qquad s,t\in \R.
\]
We study the existence of bfBm for a given pair of parameters $(\H,\K)$ and
encounter some related limiting processes.
\end{abstract}



\section{Introduction}
\setcounter{equation}{0}

Classical fractional Brownian motion (fBm) $W^{(\H)}(t), t\in\R$,
with parameter $\H\in (0,1]$, a centered Gaussian process  with
covariance
\be \label{cov_WH}
    R_{W}^{(\H)}(s,t):=
    \frac 12 \left( |s|^{2\H}+ |t|^{2\H}-|t-s|^{2\H}\right),
\ee
is so widely known and used that it needs no further recommendations.
The remarkable properties of this class of processes are described e.g.
in \cite[Section 7.2]{ST} and in \cite[Chapter 4]{EM}.

Houdr\'e and Villa \cite{HV} introduced an extension of fBm called
bifractional Brownian motion (bfBm) as a centered Gaussian process
$\b^{(\H,\K)}$ on $\R$ with covariance
\be \label{cov_BHK}
   R_\b^{(\H,\K)}(s,t):=
   2^{-\K} \left( \left(|s|^{2\H}+|t|^{2\H} \right)^{\K}-|t-s|^{2\H\K}\right).
\ee

Usual fBm shows up here when $\K=1$.

There is one more special case of bfBm directly related to the usual
fBm. Consider an anti-symmetrized version of fBm,
\[
   V^{(\H)} (t):=   W^{(\H)}(t) -   W^{(\H)}(-t), \qquad t\ge 0.
\]
It is easy to find its covariance
\[
    R_{V}^{(\H)}(s,t) = (s+t)^{2\H}-|t-s|^{2\H}, \qquad s,t\ge 0.
\]
By comparing this formula with \eqref{cov_BHK}, we see that bfBm $\b^{(1/2,\K)}$, $0<\K<2$,
consists, up to a scaling factor, of the two independent versions of
 $V^{(\K/2)}$, -- one for positive, another for negative times.\footnote{We did not find this identification with
anti-symmetrized fBm in the literature, although  the sister object,
a symmetrized fBm   $W^{(\H)}(t)+W^{(\H)}(-t)$, appears in \cite{BGT} under the name of sub-fractional Brownian motion
in connection to the limiting behavior of occupation time of particle systems.}
\bigskip

Houdr\'e and Villa motivate bfBm just by saying that "usual fBm seems to be a
valuable model for small increments [of real processes], but it appears to
be inadequate for large increments. It is thus very natural to explore the
existence of processes which keep some of the properties of fBm but also enlarge
our modelling tool kit". Marouby \cite{Mar} confirmed this deep guess by showing how
a family of bfBm's $\H=\tfrac{1}{2}$, $\K\in(0,1)$, naturally appears as a limit
in Mandelbrot micropulse model (see also \cite[Section 14]{Lif}). On the other hand,
in \cite{LST} bfBm was used for proving new probabilistic inequalities.

Initially, Houdr\'e and Villa proved the existence of bfBm on $\R$ for
\[
   \ 0<\H\le 1, \ 0< \K   <1 \, .
\]
Later on, Bardina and Es-Sebaiy \cite{BES} enlarged the zone of existence. Using an idea
of Lei and Nualart \cite{LN}, they proved that bfBm exists on $\R$ for
\[
   0<\H\le 1, 0<\K    \le \min \left\{2, \tfrac{1}{\H} \right\}.
\]
To the moment when we started this work, it was still unknown whether bfBm exists for
any other pairs $(\H,\K)$.
We  show below  in Proposition \ref{2nes} that conditions $\K\le 2$ and $\H \K \le 1$
are necessary for the existence of bfBm on $\R_+$.
\medskip

In the zone $\H>1,0<\K <\tfrac{1}{\H}$, which is most difficult for the research,
we proceed with spectral analysis and trace a new numerical bound between the zones
of existence and non-existence. We are guided by a guess of D.S.~Egorov who
conjectured that for any fixed $\H\ge 1$ there exists a positive $\bar{\K}(\H)<\H^{-1}$
such that bfBm exists for $(\H,\K)$ with any $\K<\bar{\K}(\H)$ and does not exist for
any $\K>\bar{\K}(\H)$.

\section{Existence arguments}
\setcounter{equation}{0}

For reader's convenience, we briefly recall here (and extend)
the key arguments from \cite{HV} for the case $0<\K<1$, and those from
\cite{BES, LN, Ma} for the case $1<\K\le 2$, proving the existence of bfBm.

\subsection{Case $0<K<1$}

The arguments of Houdr\'e and Villa actually have nothing to do with
fBm or bfBm, as the following statement shows.

Recall that a {\it Bernstein function} is a function $f:\R_+\to\R_+$ which
admits the following L\'evy-Khintchine representation
\be \label{Bern}
  f(\lambda) = a+b\lambda+\int_0^\infty (1-e^{-x \lambda}) \mu(dx),
\ee
where $a,b\ge 0$ are some constants and $\mu$ is a measure on $(0,\infty)$ satisfying
the integrability condition
\[
 \int \min\{x,1\} \mu(dx)<\infty.
\]
Bernstein functions, many examples and their connections to various fields of
mathematics are discussed in the monograph \cite{SSV}. Typical examples are
$\lambda\to \log(1 + \lambda)$ and
$\lambda\to \lambda^\K$ for $0<\K\le 1$. If $0<\K<1$, the representation
\eqref{Bern} takes the form
\be \label{lamK}
   \lambda^\K = \frac{\K}{\Gamma(1-\K)} \int_0^\infty (1- e^{-x\lambda} ) x^{-1-\K} dx.
\ee

\begin{prop} \label{p:exist} Let
$Y(t), t\in \R$ be a centered process with stationary increments and
finite second moments
\[
   \sigma(t)^2:=\E Y(t)^2.
\]
Then for any Bernstein function  $f(\cdot)$ there exists a process with covariance
\[
  R_{f,\sigma}(s,t):= f\left(\sigma(s)^2+\sigma(t)^2\right)-f(\sigma(s-t)^{2}),
  \qquad s,t\in \R.
\]
\end{prop}

\begin{rem} For fBm $Y=W^{(\H)}$ we have $\sigma(t)=|t|^{\H}$, thus Proposition
$\ref{p:exist}$ used with $f(\lambda)=\lambda^\K$  proves the existence of bfBm
with $0<\H,\K\le 1$.
\end{rem}

\pf
For $f(\lambda)=a+b\lambda$ we simply have
\begin{eqnarray*}
   R_{f,\sigma}(s,t) &=& b\left[
   \sigma(s)^2+\sigma(t)^2-\sigma(s-t)^{2}\right]
   \\
   &=& b \left[ \E Y(s)^2+ \E Y(t)^2 -\E (Y(t)-Y(s))^2 \right]
   \\
   &=& 2 \, b \, \E Y(s)Y(t) := 2\, b \, R_Y(s,t).
\end{eqnarray*}
Therefore, the process $\widetilde Y(t):= \sqrt{2\, b} \, Y(t)$ solves the problem.

Let now $a=b=0$. In view of the formula \eqref{Bern},
it is sufficient to find a process on $\R$ with covariance
\begin{eqnarray*}
&&
  R_x^{(\sigma)}(s,t) := \left(1- \exp(-x(\sigma(s)^{2}+\sigma(t)^{2}) \right)
  - \left(1- \exp(-x \sigma(t-s)^{2})\right)
\\
&=&  \exp(-x \sigma(t-s)^{2})   - \exp(-x(\sigma(s)^{2}+\sigma(t)^{2}))
\\
&=&   \exp(-x \sigma(s)^{2}) \exp(-x \sigma(t)^{2})
      \left[ \exp(x(\sigma(s)^{2}+\sigma(t)^{2} -\sigma(t-s)^{2}))- 1\right]
\\
&=&  \exp(-x \sigma(s)^{2}) \exp(-x \sigma(t)^{2})
     \left[ \exp(2x\, R_{Y}(s,t))- 1\right]
\\
&=&  \exp(-x \sigma(s)^{2}) \exp(-x \sigma(t)^{2})
     \sum_{m=1}^\infty \frac{(2x)^m}{m!}\, R_{Y}(s,t)^m
\end{eqnarray*}
for any $x>0$.
The latter clearly exist along with processes having covariances
$R_{Y}(\cdot,\cdot)^m $.\quad $\Box$

\subsection{Case $1<K\le 2$}

Following Lei and Nualart \cite{LN}, consider the real Gaussian process
\be \label{XKint}
  X_0^{(\K)}(t) := \int_0^\infty (1-e^{-rt}) r^{-(1+\K)/2} \W(dr),
  \qquad t \ge 0,
\ee
where $\W$ is an appropriate uncorrelated Gaussian noise. The process $X^{(\K)}$ is
well defined for $\K\in(0,2)$. By using \eqref{lamK}, and analogous formula
for $K\in(1,2)$,
\[
   \lambda^\K = \frac{\K(\K-1)}{\Gamma(2-\K)}
   \int_0^\infty (e^{-x\lambda} -1 + x\lambda ) x^{-1-\K}\, dx,
\]
it is easy to calculate the covariance
\[
   R_{X,0}^{(\K)}(s,t) :=cov\left(X_0^{(\K)}(s), X_0^{(\K)}(t)\right).
\]
We have
\[
   R_{X,0}^{(\K)}(s,t) = \begin{cases}
   \frac{\Gamma(1-\K )}{\K } \left(s^\K  + t^\K -(s+t)^\K \right), & \K \in (0,1),
   \\
   \ln s +\ln t-\ln(s+t), & \K=1,
   \\
   \frac{\Gamma(2-\K)}{\K(\K-1)} \left(-s^\K - t^\K +(s+t)^\K \right), & \K \in (1,2).
   \end{cases}
\]
Next, we rescale time by introducing a process
\be \label{XHK}
   X_0^{(\H,\K)}(t) : =  X_0^{(\K)}(|t|^{2\H}), \qquad t\in \R,
\ee
which has the covariance
\[
   R_{X,0}^{(\H,\K)}(s,t) :=cov\left(X_0^{(\H,\K)}(s), X_0^{(\H,\K)}(t)\right)
\]
given by
\begin{eqnarray*}
  &&  R_{X,0}^{(\H,\K)}(s,t)
\\
  &=& \begin{cases}
   \frac{\Gamma(1-\K )}{\K} \left(|s|^{2\H\K}+|t|^{2\H\K}-\left(|s|^{2\H}
   + |t|^{2\H}\right)^{\K}\right), & \K \in (0,1),
   \\
   2\H\left(\ln |s| + \ln |t|-\ln(|s|+|t|)\right), & \K =1,
   \\
   \frac{\Gamma(2-\K)}{\K(\K -1)} \left(-|s|^{2\H\K}-|t|^{2\H\K }
       +\left(|s|^{2\H}+|t|^{2\H}\right)^{\K}\right) , & \K \in (1,2).
   \end{cases}
\end{eqnarray*}
\medskip

If $\H\K \le 1$, consider the usual fBm $W^{(\H\K)}(t), t\in\R$, with covariance
 from \eqref{cov_WH},
\[
   R_W^{(\H\K)}(s,t)
   = \frac 12 \left( |s|^{2\H\K }+ |t|^{2\H\K }-|t-s|^{2\H\K }\right)
\]
and, for $\K \in (1,2)$, obtain bfBm just by adding  up the independent processes
\be \label{Bsum}
   \b^{(\H,\K)}(t) :=
   \sqrt{ \frac{\K (\K -1)} {2^\K \Gamma(2-\K )}}\,  X_0^{(\H,\K )}(t)
   + \sqrt{2^{1-\K }} \, W^{(\H\K )}(t), \qquad t\in \R.
\ee

For the boundary case $\K=2$ the integral representation \eqref{XKint} does not work
but we may simply define $X_0^{(2)}(t), t\ge 0$, as a degenerated random linear
process with covariance
\[
   R_{X,0}^{(2)}(s,t):= cov \left(X_0^{(2)}(s), X_0^{(2)}(t)\right)
   = 2st = \left(-s^2-t^2+(s+t)^2\right),
\]
then let again $X_0^{(\H,2)}(t):= X_0^{(2)}(|t|^{2\H})$, $t\in \R$,
as in \eqref{XHK}, and obtain
\[
    \b^{(\H,2)}(t):= 2^{-1}\, X_0^{(\H,2)}(t) +  2^{-1/2}\, W^{(2\H)}(t),
    \qquad t\in \R,
\]
whenever $0<\H\le \tfrac 12$.

In another adjacent case $\K =1$ the bfBm $\b^{(\H,1)}$ reduces to the classical
fBm $W^{(\H)}$. We hesitate to call it
a boundary case because it separates not the zones of existence and non-existence
but rather two existence zones with different properties.
\bigskip

In the zone $0<\K<1$ the representation \eqref{Bsum} does not work because
the signs in the covariance of $X^{(\K)}$ are opposite to the desired ones.
In exchange, we have a representation for fBm
\be \label{Wsum1}
    W^{(\H\K)}(t):= \sqrt{ \frac{\K } {2\Gamma(1-\K)}}\,  X_0^{(\H,\K)}(t)
              + \sqrt{2^{\K -1}}\, \b^{(\H,\K)}(t), \qquad t\in \R,
\ee
with independent processes on the right hand side.
This is equivalent to
\be \label{Bdif}
   \b^{(\H,\K)}(t) =  \sqrt{2^{1-\K }}\, W^{(\H\K )}(t)
      - \sqrt{ \frac{\K } {2^\K \Gamma(1-\K )}}\,  X_0^{(\H,\K )}(t), \qquad t\in \R.
\ee

Since $X_0^{(\K)}(\cdot)$ is a smooth process, it becomes obvious that
the local properties of $\b^{(\H,\K)}(\cdot)$ are the same as those of
fBm $W^{(\H\K)}(\cdot)$, cf. \cite{TX,W}.

We also see that if $0<\K <1$ and $\b^{(\H,\K)}(\cdot)$ exists,
then $B^{(\H\K)}(\cdot)$ exists  \cite{Ego}, which simply means $\H\K\le 1$.
In Proposition \ref{2nes} we show that $\H\K \le 1$ is necessary for
the existence of $\b^{(\H,\K)}(\cdot)$ whatever $\K$ is.

In the following we prefer to work with a modification of the processes
$X_0^{(\K)}$, $X_0^{(\H,\K)}$ having simpler covariances.
For $\K\in (0,1), \H>0$, let
\begin{eqnarray*}
  X^{(\K)}(t) &:=& \sqrt{ \frac{\K}{\Gamma(1-\K)} } \,  X_0^{(\K)}(t),
\\
   X^{(\H,\K)}(t) &:=& \sqrt{\frac{\K}{\Gamma(1-\K)} } \,  X_0^{(\H,\K)}(t).
\end{eqnarray*}
The respective covariances are
\begin{eqnarray*}
 R_{X}^{(\K)}(s,t) &=& \frac{\K} {\Gamma(1-\K)}\,  R_{X,0}^{(\K)}(s,t)
   = s^{\K }+t^{\K }-\left(s + t \right)^{\K}, \qquad s,t\ge 0,
\\
   R_{X}^{(\H,\K)}(s,t) &=& \frac{\K} {\Gamma(1-\K)}\,  R_{X,0}^{(\H,\K)}(s,t)
\\
   &=& |s|^{2\H\K }+|t|^{2\H\K }-\left(|s|^{2\H} + |t|^{2\H}\right)^{\K},
   \qquad s,t \in \R.
\end{eqnarray*}
Then \eqref{Wsum1} becomes
\be \label{Wsum2}
    \sqrt{2} \, W^{(\H\K)}(t) = X^{(\H,\K)}(t)
              + \sqrt{2^{\K}}\, \b^{(\H,\K)}(t), \qquad t\in \R,
\ee
or, in the language of covariances,
\be \label{RWsum}
  2  R_{W}^{(\H\K)}=  R_{X}^{(\H,\K)} + 2^\K  R_{\b}^{(\H,\K)}.
\ee
\medskip

Finally, notice that an extension to more general processes similar to Proposition \ref{p:exist}
is also possible for the range $1\le \K\le 2$, cf. \cite[Theorem 3.1(i)]{Ma}.

\section{Necessary conditions}
\setcounter{equation}{0}

First of all notice that we must distinguish the existence of bfBm on $\R_+$ and on
$\R$. This is very different from the case of usual fBm where condition $0<\H\le 1$
is necessary and sufficient for the existence in both cases.

\begin{prop} \label{2nes}
If bfBm exists on $\R_+$, then $\K \le 2$ and $\H\K \le 1$.
\end{prop}

\pf   
Since the covariance $R_\b^{(\H,\K)}(\cdot,\cdot)$ has the self-similarity
property
\[
   R_\b^{(\H,\K)}(cs,ct)=c^{2\H\K} R_\b^{(\H,\K)}(s,t),
\]
we may transform bfBm $\b^{(\H,\K)}$ into a stationary process by letting
\[
   U_\b^{(\H,\K)}(\tau):= e^{-\H\K \tau}\ \b^{(\H,\K)}(e^{\tau}).
\]
Stationarity of $U_\b^{(\H,\K)}$ means that its covariance function
depends only on the arguments' difference, i.e.
\[
   cov \left(U_\b^{(\H,\K)}(\tau_1),U_\b^{(\H,\K)}(\tau_2)\right)
   =: \, \tr_\b^{(\H,\K)}(\tau_2-\tau_1),
\]
where in our case
\begin{eqnarray} \nonumber
    \tr_\b^{(\H,\K)}(\tau)
    &=& e^{-\H\K\tau} R_\b^{(\H,\K)}(1, e^\tau)
\\ \label{RB}
    &=& e^{-\H\K \tau}  2^{-\K }
        \left( \left(1+ e^{2\H\tau} \right)^{\K}-| e^\tau-1|^{2\H\K}\right)
\\ \nonumber
    &=& \left( \cosh (\H\tau) \right)^{\K}-
        2^{(2\H-1)\K} \left|\sinh(\tau/2)\right|^{2\H\K}.
\end{eqnarray}
By H\"older inequality
\begin{eqnarray} \nonumber
  |\tr_\b^{(\H,\K)}(\tau)| &=& |\, cov(U_\b^{(\H,\K)}(0),U_\b^{(\H,\K)}(\tau))|
  \\ \nonumber
  &\le& \left[\E U_\b^{(\H,\K)}(0)^2\,  \E U_\b^{(\H,\K)}(\tau)^2 \right]^{1/2}
  \\ \nonumber
  &=& \left[\tr_\b^{(\H,\K)}(0)\cdot \tr_\b^{(\H,\K)}(0) \right]^{1/2}
  \\ \label{holder}
  &=&\tr_\b^{(\H,\K)}(0)=1,
\end{eqnarray}
hence, the function $|\tr_\b^{(\H,\K)}(\cdot)|$ must be bounded and must
attain its maximum at zero (this is a common property of all stationary
processes).

In our case, when $\tau\to + \infty$, in \eqref{RB} we have the expansions
\begin{eqnarray*}
  2^{-\K}  e^{-\H\K\tau}  \left(1+ e^{2\H\tau} \right)^{\K}
  &=&
  2^{-\K } e^{\H\K\tau} \left( 1+ e^{-2\H\tau} \right)^{\K}
\\
  &=&
  2^{-\K} e^{\H\K \tau} \left(1+ \K e^{-2\H\tau}(1+o(1)) \right)
\\
  &=&
  2^{-\K} e^{\H\K \tau}  +  2^{-\K} \K  e^{\H(\K -2)\tau}(1+o(1))
\end{eqnarray*}
and
\begin{eqnarray*}
  2^{-\K}  e^{-\H\K \tau} ( e^\tau-1)^{2\H\K}
  &=&
  2^{-\K}  e^{\H\K\tau} (1- e^{-\tau})^{2\H\K}
\\
  &=&
  2^{-\K}  e^{\H\K\tau} \left(1- 2\H\K e^{-\tau}(1+o(1))\right)
\\
  &=&
  2^{-\K }  e^{\H\K \tau} - \H\K  2^{1-\K } e^{(\H\K -1)\tau}(1+o(1))
\end{eqnarray*}
that yields
 \[
   \tr_\b^{(\H,\K)}(\tau)  =  2^{-\K } \K  e^{(\K -2) \H\tau}(1+o(1))
   + \H\K  2^{1-\K } e^{(\H\K -1)\tau}(1+o(1)).
 \]
Therefore, the boundedness of  $\tr(\cdot)$ implies that
both conditions  $\K\le 2$ and $\H\K\le 1$ are necessary for the existence of $\b^{(\H,\K)}$
on $\R_+$.

$\Box$  

Another argument for $\H\K\le 1$ is given in \cite[p.626]{Ma}.

\begin{prop} \label{p:covnes} The following two covariance based
necessary conditions hold.

a)  If bfBm exists on $\R$, then $\K \le \tfrac{1}{2\H-1}$.

b) If bfBm exists on $\R_+$, then $\K \le \kh(\H)$, where
\[
   \kh(\H):= \sup\left\{ \K: \sup_{\tau>0}
   \left(\left( \cosh (\H\tau) \right)^{\K}-
        2^{(2\H-1)\K} \left|\sinh(\tau/2)\right|^{2\H\K} \right)\le 1 \right\}
\]
and $\kh(\H)<\H^{-1}$ for $\H>1$.
\end{prop}
\medskip

\remark We do not have an analytic expression for the function $\kh(\cdot)$.
Some values of $\kh(\cdot)$ are given in Table \ref{tab:HKK} below.
\medskip

\pf

a) Assume that bfBm exists on $\R$. Since for its covariance we have
$R_\b^{(\H,\K)}(1,1)= R_\b^{(\H,\K)}(-1,-1)=1$,
it is true that
\[
  -1 \le  R_\b^{(\H,\K)}(1,-1)= 2^{-\K}\left[  2^\K - 2^{2\H\K}\right]= 1- 2^{(2\H-1)\K},
\]
whereas $(2\H-1)\K\le 1$.
\medskip

b) Assume that bfBm exists on $\R_+$. Then the stationary process $U_\b^{\H,\K}$ with covariance
$\tr_\b^{(\H,\K)}$ exists. Then \eqref{holder} yields
\be \label{kkk1}
   \sup_{\tau>0} \left( \left( \cosh (\H\tau) \right)^{\K}-
        2^{(2\H-1)\K} \left( \sinh(\tau/2)\right)^{2\H\K} \right)\le 1
\ee
which is equivalent to
\be \label{kkk2}
     \left( 2\cosh (\H\tau) \right)^{\K} \le
         \left( 2\sinh(\tau/2)\right)^{2\H\K} + 2^\K, \qquad \tau \ge 0.
\ee
It remains to notice that if \eqref{kkk2} holds for some value of $\K$, then it holds
for any smaller positive value of $\K$, since for any $a\in (0,1]$ we have
\[
   \left( 2\cosh (\H\tau) \right)^{a\K}
   \le \left[\left( 2\sinh(\tau/2)\right)^{2\H\K} + 2^\K\right]^a
   \le \left( 2\sinh(\tau/2)\right)^{2\H a \K} + 2^{a\K}.
\]

Finally, if $\H>1$ and $\K =\H^{-1}$, then
\[
  \lim_{\tau\to\infty}
  \left( \left( \cosh (\H\tau) \right)^{\K}-
        2^{(2\H-1)\K} \left( \sinh(\tau/2)\right)^{2\H\K} \right)
        =2^{1-\K }>1.
\]
Hence, \eqref{kkk1} fails for $\K=\H^{-1}$. Moreover, by continuity arguments, it also
fails for all $\K$ that are sufficiently close to $\H^{-1}$. It follows that $\kh(\H)<\H^{-1}$.

$\Box$  
\bigskip

The covariance criteria given in this section are quite elementary. They take into account only
2-dimensional distributions of the process. In order to get sharper results,
we need more refined spectral methods.

\section{Spectral analysis}
\setcounter{equation}{0}

\subsection{Stationary processes, covariances and spectral densities}

In addition to the self-similar processes $W^{(\H)}$, $X^{(\K)}$, $X^{(\H,\K)}$,
$\b^{(\H,\K)}$, let us introduce their stationary versions
\begin{eqnarray*}
  U_W^{(\H)}(\tau) &:=& e^{-\H\tau} W^{(\H)}(e^\tau);
  \\
  U_X^{(\K)}(\tau) &:=& e^{-\K\tau/2} X^{(\K)}(e^\tau);
   \\
  U_X^{(\H,\K)}(\tau) &:=& e^{-\H\K\tau} X^{(\H,\K)}(e^\tau);
  \\
  U_\b^{(\H,\K)}(\tau) &:=& e^{-\H\K\tau} \b^{(\H,\K)}(e^\tau).
\end{eqnarray*}

Notice that $U_W^{(\H)}$ is one of the well known versions of fractional
Ornstein--Uhlenbeck process, see e.g. \cite{BNPA,CKM}.

By the definition of $X^{(\H,\K)}$, we also have
\begin{eqnarray} \nonumber
U_X^{(\H,\K)}(\tau) &=& e^{-\H\K\tau} X^{(\K)}\left((e^{\tau})^{2\H}\right)
\\ \label{UXHK}
&=& e^{-\K(2\H\tau)/2} X^{(\K)}\left(e^{2\H\tau}\right) =  U_X^{(\K)}(2\H\tau).
\end{eqnarray}

The covariance functions corresponding to these four stationary processes are
\begin{eqnarray*}
  \tr_W^{(\H)}(\tau) &:=& \cosh(\H\tau)-2^{2\H-1}|\sinh(\tau/2)|^{2\H};
  \\
  \tr_X^{(\K)}(\tau) &:=& 2\cosh(\K\tau/2)- (2\cosh(\tau/2))^{\K};
   \\
  \tr_X^{(\H,\K)}(\tau) &:=& \tr_X^{(\K)}(2\H\tau)
  = 2\cosh(\H\K\tau)- (2\cosh(\H\tau))^{\K};
  \\
  \tr_\b^{(\H,\K)}(\tau) &:=& (\cosh(\H\tau))^\K- 2^{(2\H-1)\K} |\sinh(\tau/2)|^{2\H\K}.
\end{eqnarray*}
The basic equality \eqref{RWsum} transforms into
\be \label{WXB}
    2 \tr_W^{(\H\K)}(\tau) = \tr_X^{(\H,\K)}(\tau) + 2^\K \tr_\b^{(\H,\K)}(\tau).
\ee
\bigskip

Let us now pass to spectral representations. Recall that by inversion formula any
covariance function $\tr(\cdot)$ of a stationary process such that
$\tr\in L_1(\R)$ admits a spectral representation
\[
  \tr(\tau)=\int_{-\infty}^\infty e^{i\tau u} f(u) du, \qquad \tau\in \R,
\]
and the non-negative summable function $f(\cdot)$ is called the spectral density
of the corresponding process. We denote $f_W^{(\H)}$, $f_X^{(\K)}$,
$f_X^{(\H,\K)}$, $f_\b^{(\H,\K)}$ the spectral densities corresponding to
the respective covariance functions defined above.

Notice immediately that relation $\tr_X^{(\H,\K)}(\tau) = \tr_X^{(\K)}(2\H\tau)$ yields
\be \label{fHK}
  f_X^{(\H,\K)}(u) = \frac{1}{2\H}\, f_X^{(\K)}\left(\frac{u}{2\H}\right)\ .
\ee

\subsection{Spectral criterion for the existence of $\b^{(H,K)}$ }

\begin{prop} \label{p:compare} Let $\K\in (0,1), \H>0$. Then bfBm $\b^{(\H,\K)}$ exists on
$\R_+$ iff
\be \label{fcompare}
   f_X^{(\H,\K)} (u) \le 2\, f_W^{(\H\K)}(u), \qquad u\in\R.
\ee
\end{prop}

\pf a) Assume that \eqref{fcompare} holds. Then by \eqref{WXB}
\begin{eqnarray*}
    2^\K \tr_\b^{(\H,\K)}(\tau) &=&  2 \tr_W^{(\H\K)}(\tau) - \tr_X^{(\H,\K)}(\tau)
\\
   &=&  \int_{-\infty}^\infty e^{i\tau u}
   \left(  2\, f_W^{(\H\K)}(u) - f_X^{(\H,\K)} (u) \right) du
\\
   &:=&  \int_{-\infty}^\infty e^{i\tau u}  f(u) du,
\end{eqnarray*}
where $f(\cdot)$ is a nonnegative integrable function. It follows that a stationary process
$U_\b^{(\H,\K)}(\tau),\tau\in \R$, with covariance $\tr_\b^{(\H,\K)}(\cdot)$ exists,
and we obtain bfBm by letting
\[
  \b^{(\H,\K)}(t):=   t^{\H\K} U_\b^{(\H,\K)}(\ln t), \qquad t\ge 0.
\]
\medskip

b) Conversely, if a bfBm $\b^{\H,\K}$ exists on $\R_+$, then a stationary process
$U_\b^{(\H,\K)}(\tau)$, $\tau\in \R$, with covariance $\tr_\b^{(\H,\K)}(\cdot)$
exists. Since this function belongs to $L_1(\R)$, there
exists a non-negative spectral density $f$ such that
\[
   \tr_\b^{(\H,\K)}(\tau) = \int_{-\infty}^\infty e^{i\tau u}  f(u) du,
   \qquad \tau\in \R.
\]
By \eqref{WXB} it follows that for any $\tau\in \R$
\begin{eqnarray*}
     && \int_{-\infty}^\infty e^{i\tau u}
     \left(  2\, f_W^{(\H\K)} - f_X^{(\H,\K)} - 2^\K f \right)(u) du
 \\
    &=& 2 \tr_W^{(\H\K)}(\tau) - \tr_X^{(\H,\K)}(\tau) - 2^\K \tr_\b^{(\H,\K)}(\tau)
    =0.
\end{eqnarray*}
Since the kernel of Fourier transform is trivial, we have
\[
  2\, f_W^{(\H\K)} - f_X^{(\H,\K)} - 2^\K f  =0,
\]
Hence,
\[
  2\, f_W^{(\H\K)} - f_X^{(\H,\K)} = 2^\K f  \ge 0.
\]
$\Box$
\bigskip

The criterion of Proposition \ref{p:compare} becomes meaningful whenever
we have explicit formulae for the involved  spectral densities. They are
found below in this section.

\subsection{Spectrum associated to the Lei--Nualart process $X$}

By using the representation \eqref{XKint}, we obtain
 \begin{eqnarray*}
    U_X^{(\K)}(\tau) &=& e^{-\K\tau/2} X^{(\K)}(e^\tau)
    =   e^{-\K\tau/2} \sqrt{\frac{\K}{\Gamma(1-\K)}}\,  X_0^{(\K)}(e^\tau)
 \\
    &=& \sqrt{\frac{\K}{\Gamma(1-\K)}}\,
    \int_0^\infty e^{-\K\tau/2} (1-e^{-re^\tau}) r^{-(1+\K)/2} \W(dr).
 \end{eqnarray*}
 It follows that
 \begin{eqnarray*}
    \tr_X^{(\K)}(\tau) &=& cov(U_X^{(\K)}(\tau),U_X^{(\K)}(0))
  \\
    &=& \frac{\K}{\Gamma(1-\K)}
    \int_0^\infty e^{-\K\tau/2} (1-e^{-re^\tau}) (1-e^{-r}) r^{-(1+\K)} dr
  \\
   &=& \frac{\K}{\Gamma(1-\K)}
   \int_{-\infty}^\infty e^{-\K\tau/2} (1-e^{-e^{(v+\tau)}}) (1-e^{-e^v}) e^{-\K v} dv
  \\
     &=& \frac{\K}{\Gamma(1-\K)}
     \int_{-\infty}^\infty e^{-\K(v+\tau)/2} (1-e^{-e^{(v+\tau)}})
     e^{-\K v/2} (1-e^{-e^v}) dv
  \\
    &=& \frac{\K}{\Gamma(1-\K)} \int_{-\infty}^\infty g(v+\tau)g(v) dv,
\end{eqnarray*}
 where $g(v):=e^{-\K v/2} (1-e^{-e^v})$. By applying Fourier transform,
 we obtain
\begin{eqnarray*}
   \tr _X^{(\K)}(\tau) &=& \frac{\K}{\Gamma(1-\K)}  \int_{-\infty}^\infty
   e^{i\tau u} \, \widehat{g}(u) \, \overline{\widehat{g}(u)}\, du
   \\
   &=&  \frac{\K}{\Gamma(1-\K)} \int_{-\infty}^\infty e^{i\tau u}  |\widehat{g}(u)|^2  du.
\end{eqnarray*}
It follows that
\[
   f_X^{(\K)}(u)=  \frac{\K}{\Gamma(1-\K)}\,  |\widehat{g}(u)|^2.
\]
Now we find $\widehat{g}(u)$. By definition
\begin{eqnarray*}
\widehat{g}(u) &=& \frac{1}{\sqrt{2\pi}} \int_{-\infty}^\infty e^{-iu v} g(v)  dv
\\
&=& \frac{1}{\sqrt{2\pi}} \int_{-\infty}^\infty e^{-iu v} e^{-\K v/2}(1-e^{-e^v})dv
\\
&=& \frac{1}{\sqrt{2\pi}} \int_{0}^\infty r^{-iu-\K/2 -1}(1-e^{-r})dr
\\
&=& - \frac{1}{\sqrt{2\pi}(-i u-\K/2)} \int_{0}^\infty r^{-i u-\K/2}e^{-r}dr
\\
&=& - \frac{\Gamma(-i u-\K/2+1)}{\sqrt{2\pi}(-i u-\K/2)}
= - \frac{\Gamma(-i u-\K/2)}{\sqrt{2\pi}} \, .
\end{eqnarray*}
We conclude that
\[
   f_X^{(\K)}(u) =  \frac{\K}{\Gamma(1-\K)}\, \frac{|\Gamma(-iu-\K/2)|^2}{2\pi}\, .
\]
Finally,  equation \eqref{fHK} yields
\[
   f_X^{(\H,\K)}(u) = \frac{1}{2\H}\,
   \frac{\K}{\Gamma(1-\K)}\, \frac{|\Gamma( \frac{-iu}{2\H}-\K/2)|^2}{2\pi}\, .
\]

\subsection{Spectrum of the fractional Ornstein--Uhlenbeck process}

The layout of calculation is very much the same as for the spectrum of
Lei--Nualart process. Recall that fractional Brownian motion, as a process with
stationary increments, admits, for $\H\in(0,1)$, a spectral representation
\[
  W^{(\H)}(t)
  =\int_{-\infty}^{\infty} \frac{y_\H\left(e^{itr}-1 \right)}{|r|^{\H+1/2}}\ \W(dr),
\]
where
\[
   y_\H^2 = \frac{\Gamma(2\H+1)\sin(\pi \H)}{2\pi}\, .
\]
Therefore,
\[
   U_W^{(\H)}(\tau) = e^{-\H\tau}\,  W^{(\H)}(e^\tau)
  =\int_{-\infty}^{\infty} e^{-\H\tau}\,
   \frac{y_\H \left(e^{ie^\tau r}-1 \right)}{|r|^{\H+1/2}}\, \W(dr),
\]
and
\begin{eqnarray*}
   \tr_W^{(\H)}(\tau)
  &=&  y_\H^2 \int_{-\infty}^{\infty} e^{-\H\tau}\,
   \frac{ \left(e^{ie^\tau r}-1 \right)}{|r|^{\H+1/2}}
   \overline{\frac{ \left(e^{ir}-1 \right)}{|r|^{\H+1/2}}}\ dr
 \\
   &=&
   2\, y_\H^2 \ Re \int_{0}^{\infty} e^{-\H\tau}\,
   \frac{ \left(e^{ie^\tau r}-1 \right)}{r^{\H+1/2}}
   \overline{\frac{ \left(e^{ir}-1 \right)}{r^{\H+1/2}}}\ dr
\\
   &=&
   2\, y_\H^2 \ Re \int_{-\infty}^{\infty} e^{-\H(v+\tau)}\,
   \left(e^{ie^{v+\tau}}-1 \right)
   \overline{ e^{-\H v} \left(e^{ie^v}-1 \right)}\  dv
\\
   &=&  2\, y_\H^2\ Re \int_{-\infty}^{\infty} g(v+\tau)
  \, \overline{g(v)} \, dv,
\end{eqnarray*}
where $g(v):= e^{-\H v} \left( e^{ie^v}-1 \right)$.

By applying Fourier transform,
 we obtain
\begin{eqnarray*}
   \tr _W^{(\H)}(\tau) &=& 2\, y_\H^2\ Re   \int_{-\infty}^\infty
   e^{i\tau u} \, \widehat{g}(u) \, \overline{\widehat{g}(u)}\, du
   \\
   &=&   2\, y_\H^2\ Re \int_{-\infty}^\infty e^{i\tau u}  |\widehat{g}(u)|^2 \ du
   \\
   &=&   2\, y_\H^2\  \int_{-\infty}^\infty \cos(\tau u)  |\widehat{g}(u)|^2 \ du
   \\
   &=&   2\, y_\H^2\  \int_{-\infty}^\infty \cos(\tau u)
         \frac{ |\widehat{g}(u)|^2+ |\widehat{g}(-u)|^2}{2} \ du
   \\
   &=&   2\, y_\H^2\  \int_{-\infty}^\infty e^{i\tau u}
         \frac{ |\widehat{g}(u)|^2+ |\widehat{g}(-u)|^2}{2} \ du.
\end{eqnarray*}
It follows that
\[
   f_W^{(\H)}(u)= y_\H^2\  \left(|\widehat{g}(u)|^2+ |\widehat{g}(-u)|^2\right).
\]
Now we find $\widehat{g}(u)$. By definition
\begin{eqnarray*}
\widehat{g}(u) &=& \frac{1}{\sqrt{2\pi}} \int_{-\infty}^\infty e^{-iu v} g(v)  dv
\\
&=& \frac{1}{\sqrt{2\pi}}
\int_{-\infty}^\infty e^{-(iu+\H) v}\left(e^{ie^{v}}-1\right)\, dv
\\
&:=& \frac{1}{\sqrt{2\pi}}
\int_{-\infty}^\infty e^{- z v}\left(e^{ie^{v}}-1\right)\, dv
\\
&=& \frac{i}{\sqrt{2\pi} z} \int_{-\infty}^\infty e^{-z v}\, e^{v} e^{ie^{v}} \ dv
\\
&=& \frac{i}{\sqrt{2\pi}z} \int_{0}^\infty r^{-z}\,  e^{ir} \ dr
\\
&=& \frac{i}{\sqrt{2\pi}z} \int_{0}^\infty r^{-z}\,  \left(\cos r+ i\sin r\right) \ dr
\\
&=& \frac{i}{\sqrt{2\pi}z}
\left( \frac{\pi}{2\Gamma(z)\cos(\pi z/2)} +  \frac{i\pi}{2\Gamma(z)\sin(\pi z/2)}
\right)
\\
&=& \frac{-\sqrt{\pi/2} }{z\Gamma(z)} \ \frac{\cos(\pi z/2)- i\sin(\pi z/2)}{\sin (\pi z)}
\\
&=& \frac{-\sqrt{\pi/2} }{z\Gamma(z)} \ \frac{e^{-i\pi z/2}}{\sin (\pi z)} \, ,
\end{eqnarray*}
where $z= \H+iu$.
Notice that the integrals that appear after the integration by parts must be understood
as the main value integrals (they are not absolutely converging).

We infer that
\[
  |\widehat{g}(u)|^2
  =  \frac{\pi/2 }{|z|^2|\Gamma(z)|^2} \ \frac{e^{\pi u}}{|\sin (\pi z)|^2}\, .
\]
By using trigonometric formulae
\begin{eqnarray*}
  \sin (\pi z) &=& \sin (\pi \H+i\pi u)
  =   \sin(\pi \H) \cos(i\pi u) +    \cos(\pi \H) \sin(i\pi u)
\\
 &=&  \sin(\pi \H) \cosh(\pi u) +  i \cos(\pi \H) \sinh(\pi u),
\end{eqnarray*}
\[
  |\sin (\pi z)|^2 = \sin^2(\pi \H) \cosh^2(\pi u) +  \cos^2(\pi \H) \sinh^2(\pi u),
\]
we may conclude that
\begin{eqnarray*}
 f_W^{(\H)}(u) &=& y_\H^2\   \frac{\pi/2 }{|z|^2|\Gamma(z)|^2} \
 \frac{e^{\pi u}+ e^{-\pi u}}{\sin^2(\pi \H) \cosh^2(\pi u) +  \cos^2(\pi \H) \sinh^2(\pi u)}
 \\
 &=&    \frac{\Gamma(2\H+1)\sin(\pi \H)}{2 (u^2+\H^2) |\Gamma(\H+iu)|^2} \
 \frac{\cosh(\pi u)}{\sin^2(\pi \H) \cosh^2(\pi u) +  \cos^2(\pi \H) \sinh^2(\pi u)}.
\end{eqnarray*}

An interesting special case is $\H=1/2$ where, using a formula 6.1.30 from
\cite{AS},
\[
   |\Gamma(1/2+iu)|^2= \frac{\pi}{\cosh(\pi u)} \, ,
\]
we get
\[
   f_W^{(1/2)}(u) =  \frac{2}{\pi(4u^2+1)}\ ,
\]
in accordance with the classical Ornstein-Uhlenbeck covariance function
\[ \tr_W^{(1/2)}(\tau)= e^{-|\tau|/2}.
\]
\medskip

There is an alternative approach to the computation of the spectral density
$f_W^{(\H)}$, $\H\not=1/2$, due to Barndorff-Nielsen and Perez-Abreu, \cite{BNPA}.
Writing
\[
  \tr_{W}^{(\H)}(\tau) = \frac{e^{\H\tau}}{2}
  \left(1+e^{-2\H\tau} - (1-e^{-\tau})^{2\H}\right),
  \qquad \tau>0,
\]
and using Taylor expansion
\[
 (1-v)^{2\H} = 1+\sum_{k=1}^\infty (-1)^k C_{k,\H} \, v^k, \qquad
 C_{k,\H} :=\frac{\Gamma(2\H+1)}{\Gamma(k+1)\Gamma(2\H-k+1)}\, ,
\]
one obtains
\[
  \tr_{W}^{(\H)}(\tau)
  = \frac{e^{-\H|\tau|}}{2} +  \sum_{k=1}^\infty (-1)^{k+1} \frac{C_{k,\H}}{2}
  \, e^{-(k-\H)|\tau|}    \, , \qquad \tau\in\R,
\]
whereas
\[
  f_W^{(\H)}(u)
  = \sum_{k=0}^\infty (-1)^{k+1} C_{k,\H} \frac{k-\H}{2\pi(u^2+(k-\H)^2)}\, .
\]

\section{Computation}
\setcounter{equation}{0}
\setcounter{table}{0}

According to Proposition \ref{p:covnes}, in order to check the existence of bfBm on $\R_+$,
one must check for each pair $(\H,\K)$ with $\H>1$, $0<\K<\H^{-1}$, whether for all $u\ge 0$,
\begin{eqnarray} \nonumber
&&  f_X^{(\H,\K)}(u)
= \frac{1}{2\H}\, \frac{\K}{\Gamma(1-\K)}\, \frac{|\Gamma( \frac{-iu}{2\H}-\K/2)|^2}{2\pi}
\\ \label{comppare}
&\le& 2 f_W^{(\H\K)}(u)
\\
&=&  \nonumber
 \frac{\Gamma(2\H\K+1)\sin(\pi \H\K)}{ (u^2+(\H\K)^2) |\Gamma(\H\K+iu)|^2} \
 \frac{\cosh(\pi u)}{\sin^2(\pi \H\K) \cosh^2(\pi u) +  \cos^2(\pi \H\K) \sinh^2(\pi u)}.
\end{eqnarray}

Our numerical computations show that for every fixed $\H>1$  there is a positive bound $\bar \K$ such that
condition  \eqref{comppare} holds for all $\K<\bar \K$ and does not hold for all $\K>\bar \K$.

More precisely, we went through the range $\H\in (1,100]$ with a step $0.01$. For every $\H$ we went through
the range $\K\in(0,\H^{-1})$ with the same step $0.01$ and checked inequality \eqref{comppare} for $u\in [0,50]$
with the step $0.01$. Computations indicate that larger values of $u$ are by far irrelevant. They also show the existence
of the boundary $\bar{\K} :=\bar{\K}(\H)$, separating the existence and the non-existence zones as said above. Furthermore,
we repeated the procedure with the smaller step $0.001$ in $\K$ near the boundary value, in order to get sharper
values of function $\bar{\K}$.

Few words about the boundary case: since the expression for covariance function of bfBm is continuous in $H$ and $K$,
it is clear that the property of its non-negative definiteness is conserved when we let the parameters go to some limits.
In other words, the domain of existence of bfBm is closed on the plane $(H,K)$. In particular, the bfBm
with boundary parameters $(H,\overline{K}(H))$ should exist.
\medskip

A sample of values of  $\bar{\K}(\H)$ is given in Table \ref{tab:HKK} along with its upper bound $\kh(\H)$.
The complete table of values of functions $\bar{\K}(\H)$ and $\kh(\H)$ may be found in \cite{LV_arxiv}. The
bound  $\kh(\H)$ is rather sharp, especially for large values of $\H$.
\medskip

\begin{table}
{\small
\noindent\begin{tabular}{|c|c|c||c|c|c||c|c|c|}
  \hline
  $\H$ &  $\bar{\K}(\H)$ & $\kh(\H)$ & $\H$ &  $\bar{\K}(\H)$ & $\kh(\H)$ &$\H$ &  $\bar{\K}(\H)$ & $\kh(\H)$   \\
  \hline
  1.01  & 0.988 & 0.988 & 2 & 0.422 & 0.440 & 6 & 0.117 & 0.123 \\
  \hline
   1.1  & 0.887 & 0.894 & 2.5 & 0.321 & 0.338 & 7 & 0.099 & 0.104 \\
  \hline
   1.2  & 0.794 & 0.807 & 3 & 0.260 & 0.273 & 10 & 0.067 & 0.070 \\
  \hline
   1.3  & 0.718 & 0.734 & 3.5 & 0.217 & 0.228 & 20 & 0.032 & 0.033 \\
  \hline
   1.5 & 0.603 & 0.619 & 4 & 0.185 & 0.196 & 60 & 0.010 & 0.010 \\
  \hline
   1.7 & 0.519 & 0.533 & 5 & 0.144 & 0.152 & 100 & 0.006 & 0.006 \\
  \hline
\end{tabular}

}
\caption{Existence boundary $\bar{\K}(\H)$ and its covariance upper bound $\kh(\H)$}\label{tab:HKK}
\end{table}
\medskip

\begin{center}
\begin{figure} [ht]

\begin{center}
\includegraphics[height=2.5in,width=4in]{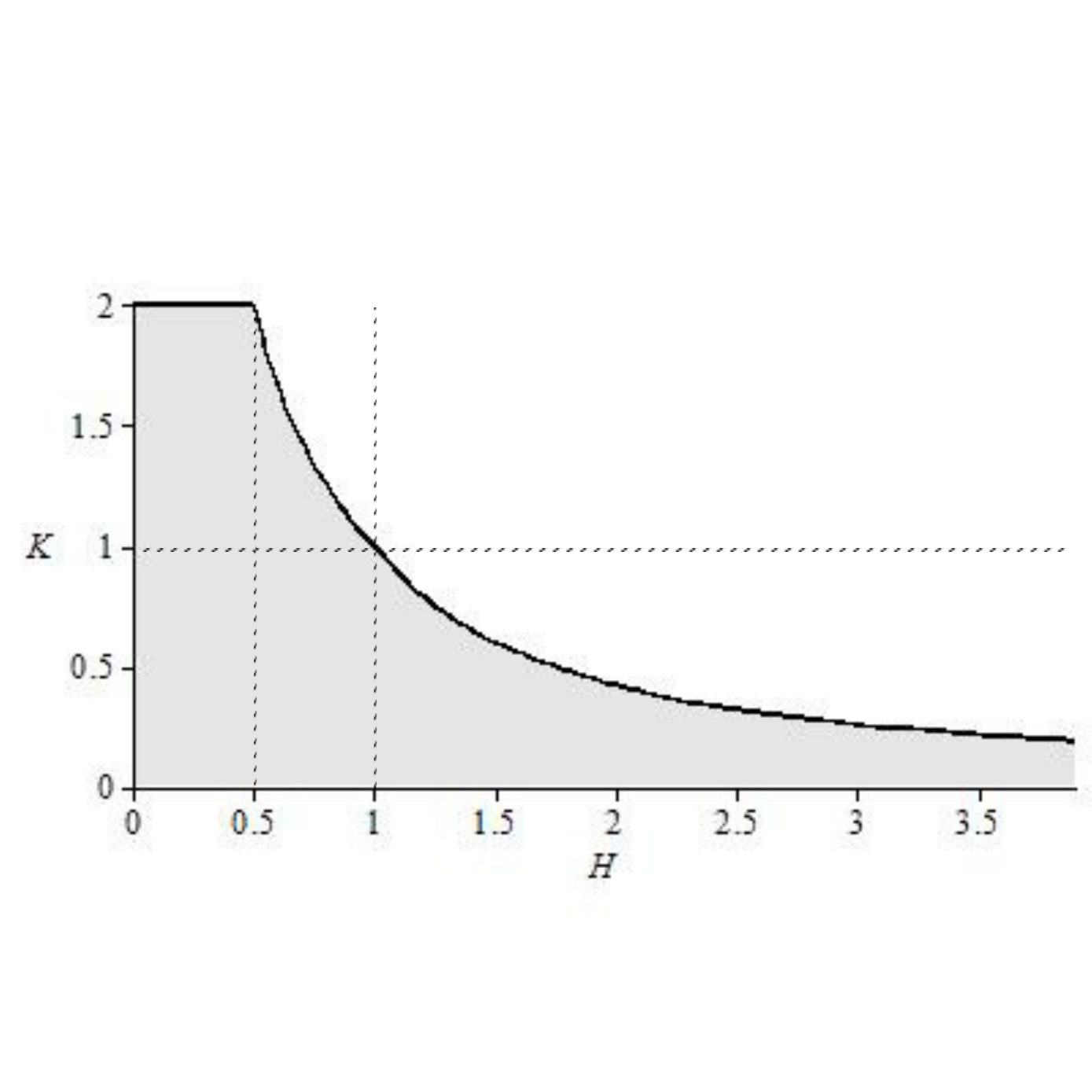}
\end{center}

\caption{Global zone of existence for bfBm $\{\b^{(\H,\K)}(t), t\ge 0\}$.}
\label{fig:exist}
\end{figure}
\end{center}

The resulting global zone of existence for bfBm $\{\b^{(\H,\K)}(t), t\ge 0\}$,
is represented at Figure \ref{fig:exist}.

\section{Some boundary cases}
\setcounter{equation}{0}

\subsection{A limiting process for $K=0$}

Let us consider a limiting behavior of the bfBm covariance function when $\H>0$ is
fixed  and $\K\to 0$. For $\tau>0$ we have
\begin{eqnarray*}
&& \K^{-1} \tr^{(\H,\K)}_\b (\tau)
\\
&=& \K^{-1} \left( \exp\left(\K \ln(\cosh(\H\tau))\right)
    -\exp\left(\K((2\H-1)\ln 2 +\ln\sinh(\tau/2)) \right)\right)
\\
&\to&  \ln(\cosh(\H\tau)) -  (2\H-1)\ln 2 - \ln\sinh(\tau/2)
\\
&=&  \ln(2 \cosh(\H\tau)) -  2\H \ln(2\sinh(\tau/2))
\\
&=&  \left[ \ln(2 \cosh(\H\tau)) -\H\tau \right]
     + 2\H \left[ \tau/2 - \ln(2\sinh(\tau/2)) \right]
     := R_1(\tau)+ R_2(\tau).
\end{eqnarray*}
We want to find the spectrum corresponding to this limiting
covariance\footnote{Notice however a logarithmic explosion of the term $R_2$
at zero. This means that the limiting process is not a usual process defined
pointwise but a generalized one. This feature may not be repaired by time scaling.}.
Let
\[
  f_j(u):=\frac{1}{2\pi} \int_{-\infty}^{\infty} e^{-i u \tau} R_j(\tau) \, d\tau,
  \qquad j=1,2,
\]
denote the corresponding spectral densities. In order to find the densities $f_j$,
we use the classical relation between the differentiation and the Fourier transform,
\[
   f_j(u)  = \frac{\widehat{R_j'}(u)}{\sqrt{2\pi}(iu)}\, , \qquad j=1,2.
\]
Since
\[
   R'_1(\tau)= \frac{\H(e^{\H\tau}-e^{-\H\tau})}{e^{\H\tau}+e^{-\H\tau}} -\H\sgn(\H\tau)
   = \H(\tanh(\H\tau) -\sgn(\H\tau)),
\]
we have
\[
  \widehat{R'_1}(u) =\widehat{\tanh-\sgn}(u/\H).
\]
Furthermore, since
\[
  \widehat{\tanh-\sgn}(u) = - i \sqrt{\pi/2}\, (\sinh(\pi u/2))^{-1}
  + i \sqrt{2/\pi}\,  u^{-1},
\]
we obtain
\begin{eqnarray} \nonumber
  f_1(u) &=&  \frac{\widehat{\tanh-\sgn}(u/H)}{\sqrt{2\pi}(iu)}
  \\  \label{f1H0}
  &=& -\left(2 u \sinh(\pi u/2\H)  \right)^{-1} + \frac{\H}{\pi u^2}
  := H^{-1} \phi(u/\H),
\end{eqnarray}
where
\[
  \phi(u) = \pi^{-1} u^{-2} \left[
  1 - \frac{\pi u}{2\sinh(\pi u/2)}\right].
\]
This is a nice function with finite limit at zero and quadratic decay at infinity.
\medskip

Next, easy calculation shows that for $\tau>0$
\begin{eqnarray*}
   R'_2(\tau)&=& \H \left( 1- \coth(\tau/2) \right)
\\
 &=& \frac{-2\H e^{-\tau}}{1-e^{-\tau}} = - 2\H \sum_{n=1}^\infty e^{-n \tau}.
\end{eqnarray*}
Hence,
\[
  R'_2(\tau) = - 2\H \sgn(\tau) \sum_{n=1}^\infty e^{- n|\tau| n}, \qquad \tau\in \R,
\]
whereas
\[
  \widehat{R'_2}(u) = \frac{4\H i u}{\sqrt{2\pi}} \sum_{n=1}^\infty (u^2+n^2)^{-1}
\]
and
\[
  f_2(u) =   \frac{\widehat{R'_2}(u)}{\sqrt{2\pi}(iu)}
  = \frac{2\H}{\pi} \sum_{n=1}^\infty (u^2+n^2)^{-1}.
\]
The spectral density $f_2(\cdot)$ is locally nice but it decays like
$|u|^{-1}$ at infinity. Therefore, it is not integrable and corresponds to
a {\it generalized} Gaussian process.
\medskip

By summing up, we obtain the spectral density
\[
  f(u) = f_1(u)+f_2(u)
  =
  -\left(2 u \sinh(\pi u/2\H)  \right)^{-1} + \frac{\H}{\pi u^2}
  +  \frac{2\H}{\pi} \sum_{n=1}^\infty (u^2+n^2)^{-1}.
\]
\bigskip

Presence of hyperbolic functions in the computations suggests that there
should be some relation of the introduced objects to hyperbolic geometry.
This is indeed the case. Cohen and Lifshits studied in \cite{CL}
many random fields and processes on the hyperbolic space. In particular,
they introduced so called quadratic field playing important role in a hyperbolic
version of spectral representations. As shown in \cite[Section 10.1]{CL},
being restricted on a geodesic line of the hyperbolic plane, quadratic field
generates a centered Gaussian process with stationary increments $Z(\tau)$,
$\tau\in\R$, with the structure function
\[
   \E \, Z(\tau)^2=2 \ln\cosh(\tau/2), \qquad tau \in \R.
\]
Notice by the way that the derivative $Z'(\cdot)$ is a stationary process
with the spectral density
\[
   f_Z(u)=\frac {u}{2\sinh(\pi u)}.
\]
A similar expression already appeared in \eqref{f1H0}.

Let us fix $H=1$ and denote $\b_1^{(1,0)}(\tau)$ a stationary Gaussian
process with spectral density \eqref{f1H0}.

Then straightforward calculations show that independent copies of $Z(\cdot)$
and $\b_1^{(1,0)}(\cdot)$ are connected by
\[
   Z(2\tau) +\b_1^{(1,0)}(\tau)-\b_1^{(1,0)}(0) = \sqrt{2}\, W(\tau),
\]
where $W=W^{(1/2)}$ is a Wiener process.

\subsection{Case $H=1$: integral representation}

Very few is known about white noise integral representations of bfBm
(for other processes, see e.g. \cite[Section 7.3]{Lif}).
We present here one for the boundary case $\H=1$, $0<\K<1$. By using
\eqref{lamK}, we have
\begin{eqnarray*}
   R_\b^{(1,\K)}(s,t)
   &=&
   2^{-\K} \left( \left(s^{2}+t^{2} \right)^{\K}-(t-s)^{2\K}\right)
\\
   &=& \frac{\K}{2^{\K} \Gamma(1-\K)}  \int_0^\infty
   \left(   e^{-x(t-s)^2}  -  e^{-x(s^{2}+t^{2})}    \right)   x^{-1-\K} dx
\\
   &=& \frac{\K}{2^{\K} \Gamma(1-\K)}  \int_0^\infty
   \left(  e^{-x(s^{2}+t^{2})} \left(  e^{2xst}-1  \right)\right)x^{-1-\K} dx
\\
   &=& \frac{\K}{2^{\K} \Gamma(1-\K)}  \sum_{n=1}^\infty  \int_0^\infty
   \left(  e^{-x(s^{2}+t^{2})}  \frac{(2xst)^n}{n!}\right)x^{-1-\K} dx.
\end{eqnarray*}
It follows that
\[
  \b^{(1,\K)}(s)=  \left(\frac{\K}{2^{\K} \Gamma(1-\K)}\right)^{1/2}
  \sum_{n=1}^\infty  s^n \, \left(\frac{2^n}{n!}\right)^{1/2}  \int_0^\infty
  e^{-x s^{2}} \, x^{(n-1-\K)/2} \W_n(dx),
\]
where $\W_n$ are independent Gaussian white noises on $\R_+$ controlled by
Lebesgue measure.

\section{A posterior discussion}

One of the referees made some remarks  related to a more general context than
the particular results of this note. We also believe that a wider discussion
might be interesting to the reader.

The referee states ``Since the introduction of bfBm there is a real problem of motivation.
Apart from self-similarity, this process enjoys no inherent property, and there
are lots of processes with the same qualitative behaviors: a lot of Bernstein functions
and functions $\sigma$ can be used and, with the help of of Proposition 2.1, new
processes may be introduced''. We basically agree with that, although we find that
self-similarity is quite a strong additional feature for picking bfBm from the crowd
of possible generalizations of fBm. In our opinion, the decisive argument for motivation
of the interest would be finding some natural models converging to the studied process.
So far, only the case $\H=\tfrac{1}{2}$ of bfBm appeared in this setting but, as we noticed
(may be for the first time?), this case is directly related to fBm itself and, therefore,
may not advocate for the entire family of bfBm's.

Back to bfBm, the referee remarks correctly that our existence results are not the same for $\R$ and
$\R_+$. This is of course the weakness of our spectral approach (that is focused only on $\R_+$) and
should be considered as a source of open problems. Moreover, considering multivariate case only makes
sense when the existence of bfBm on the whole $\R$ will be understood. Solving these problems apparently
goes far beyond the means of our note.

\bigskip

{\bf Acknowledgement.}
We are very grateful to both anonymous referees for the careful reading of our note and for
their insightful remarks.

\end{document}